\documentclass[10pt]{amsart}

\usepackage{amssymb,amsmath,amsthm}
\usepackage{graphicx}

\theoremstyle{plain}

\usepackage{graphicx,tikz}
\newtheorem{theorem}{Theorem}

\newtheorem{corollary}{Corollary}

\theoremstyle{definition}

\theoremstyle{remark}

\newcommand{\vol}{\operatorname{vol}}

\begin{document}

\title[]{Oscillatory functions vanish on a large set}

\author[]{Stefan Steinerberger}
\address{Department of Mathematics, Yale University}
\email{stefan.steinerberger@yale.edu}
\thanks{This work is supported by the NSF (DMS-1763179) and the Alfred P. Sloan Foundation.}

\begin{abstract} Let $(M,g)$ be a $n-$dimensional, compact Riemannian manifold. We define the frequency scale $\lambda$ of a function  $f \in C^{0}(M)$ as the largest number such that $\left\langle f, \phi_k \right\rangle =0$ for all Laplacian eigenfunctions with eigenvalue $\lambda_k \leq \lambda$. If $\lambda$ is large, then the function $f$ has to vanish on a large set
$$ \mathcal{H}^{n-1} \left\{x:f(x) =0\right\} \gtrsim_{}  \left( \frac{ \|f\|_{L^1}}{\|f\|_{L^{\infty}}} \right)^{2 - \frac{1}{n}}    \frac{  \sqrt{\lambda}}{(\log{\lambda})^{n/2}}.$$
Trigonometric functions on the flat torus $\mathbb{T}^d$ show that the result is sharp up to a logarithm if $\|f\|_{L^1} \sim \|f\|_{L^{\infty}}$. We also obtain a stronger result conditioned on the geometric regularity of $\left\{x:f(x) = 0\right\}$. 
This may be understood as a very general higher-dimensional extension of the Sturm oscillation theorem. 
\end{abstract}

\maketitle

\vspace{-10pt}

\section{Introduction}
\subsection{Introduction.}
 The classical Sturm oscillation theorem states that if we consider the operator $-\Delta + V$ on $[a,b]$ with
Dirichlet boundary conditions, then the $n-$th eigenfunction has exactly $n-1$ roots and a small set of zeroes is only achieved by low frequencies.
This statement is now well embedded in a much bigger context (see e.g. \cite{anton} on Sturm-Liouville theory).
One interesting result in this spirit, for $V=0$, is the following: if $f \in C^{0}(\mathbb{T})$ is
orthogonal to $\left\{1, \sin{x}, \cos{x}, \sin{2x}, \dots, \cos{n x}\right\}$, then $f$ has at least $2n+2$ roots
(for the proof: identify $\mathbb{T} = \partial \mathbb{D}$, consider the Poisson extension, conclude that the nodal set in the origin has at least $n+1$ lines and use
the maximum principle to see that these lines cannot intersect; more quantitative versions can be found in \cite{stein2, stein3}). These results apply to the eigenfunctions of general Sturm-Liouville operators, however, there is no known form of Sturm-Liouville theory in higher extension.
\begin{quote}
Even the Sturm theory is missing in higher dimensions. This is an interesting phenomenon. All the attempts I know to extend Sturm theory to higher dimensions failed. (V. I. Arnold, Third Lecture at 1997 Conference in Honor of his 60th birthday \cite{arnold})
\end{quote}

Any generalization of Sturm Theory to higher dimensions or compact Riemannian manifolds $(M,g)$ is bound to be highly nontrivial: the natural analogue of
the trigonometric functions are the Laplacian eigenfunctions. However, already the simple question of how much
a single eigenfunction has to vanish is the subject of a long-standing conjecture of S.-T. Yau \cite{yau} (see also \cite{nadirashvili})
$$ \mathcal{H}^{n-1} \left\{x \in M: \phi_k(x) = 0\right\} \sim \sqrt{\lambda_k}.$$
The question has inspired a large body of work \cite{brun, colding, don, don2, lin0, lin1, hardt, hezari, hezari2, lin2, log2, mal, mangoubi, mangoubi2, nadirashvili, sog1, sog2, stein1} and the lower bound was recently established by Logunov \cite{logunov}. On the one-dimensional torus $\mathbb{T}$, the question is whether $\sin{kx}$ has $\sim k$ roots while the Sturm-Hurwitz Theorem mentioned above implies that \textit{any} linear
combination of eigenfunctions with eigenvalue $\lambda_k \gtrsim \lambda$ has $\gtrsim \sqrt{\lambda}$ roots 
$$ \# \left\{x \in \mathbb{T}:  \sum_{k^2 \gtrsim \lambda}{a_k \sin{k x}} = 0\right\} \gtrsim \sqrt{\lambda}.$$
We were interested in whether this is a more general phenomenon: is it possible to ask the question not just for a single
eigenfunction but for the entire orthogonal complement of the first few eigenfunctions? Or, phrased differently, is it
true that functions orthogonal to the first few
Laplacian eigenfunctions vanish on a large set? This is indeed the case and we present a first quantitative result in that direction.

 \subsection{The Main Result.} The purpose of this paper is to show that results of this type are indeed possible. 
 We will work on a smooth, compact, $n-$dimensional Riemannian manifold $(M,g)$ (with or without boundary) and continuous functions $f \in C^0(M)$. 
Such functions admit
 a decomposition into eigenfunctions of the Laplacian $-\Delta_g \phi_k = \lambda_k \phi_k$ (which we assume to be $L^2-$normalized $\|\phi_k\|_{L^2}=1$)
 $$ f = \sum_{k \geq 0}{ \left\langle f, \phi_k \right\rangle \phi_k}.$$
  We define the frequency scale $\lambda \in \mathbb{R}_{\geq 0}$ of a function $f$ as the largest real number such
that $f$ is \textit{almost} orthogonal to all eigenfunctions whose eigenvalue is $\lambda_k \leq \lambda$. One
could demand complete orthogonality but the proof has a little bit of wiggle room and it suffices to assume
$$ \left(\sum_{\lambda_k \leq \lambda}{ |\left\langle f, \phi_k \right\rangle|^2}\right)^{1/2} \leq c \| f\|_{L^1},$$
where $0<c \ll 1$ is a constant only depending on the dimension.
We emphasize that the norm on the right-hand side is $L^1$ and not, as one might assume from Hilbert space
geometry, $L^2$. It seems conceivable that a version of the statement is true with $L^2$ instead of $L^1$ but our
proof does not show that.
Summarizing, $\lambda$ is defined so that $f$ is almost orthogonal to all Laplacian eigenfunctions with eigenvalue $\lambda_k \leq \lambda$ and lives in the high-frequency spectrum.

\begin{center}
\begin{figure}[h!]
\begin{tikzpicture}[scale=4.5]
\draw [ultra thick] (0,0) -- (1,0) -- (1,1) -- (0,1) -- (0,0);
\draw [ultra thick] (0, 0.3) to[out=0, in=120] (0.4, 0);
\draw [ultra thick] (0.4, 1) to[out=330, in=180] (1, 0.3);
\draw [ultra thick] (0,1) -- (1,0);
\draw [ultra thick] (0,0.6) -- (0.6,0);
\draw [ultra thick] (0,0) -- (1,1);
\node at (0.2, 0.1) {$+$};
\node at (0.1, 0.34) {$+$};
\node at (0.05, 0.15) {$-$};
\node at (0.6, 0.2) {$+$};
\node at (0.9, 0.7) {$+$};
\node at (0.7, 0.9) {$-$};
\node at (0.2, 0.6) {$-$};
\node at (0.5, 0.8) {$+$};
\node at (0.6, 0.5) {$-$};
\node at (0.4, 0.1) {$-$};
\end{tikzpicture}
\caption{A function $f$ with this sign pattern on $[0,1]^2$ is probably not orthogonal to the first 1000 eigenfunctions.}
\end{figure}
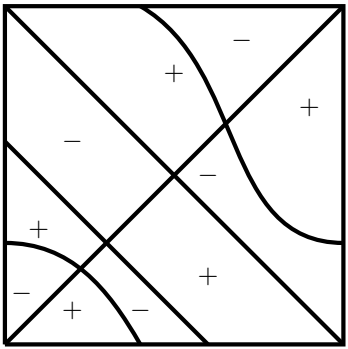
\end{center}

Motivated by classical heuristics, one would expect that such a sum of high-frequency components should not be able
to avoid vanishing on a large set. This is somewhat dual to work of Donnelly \cite{donn} showing that a sum of finitely
many eigenfunctions cannot vanish subtantially more than the eigenfunction associated to the largest eigenvalue (see
also Jerison \& Lebeau \cite{jer}).\\

The reverse statement (see Fig. 1) is also interesting: a function $f$ vanishing in a small set is probably not orthogonal to the first $k$ eigenfunctions with $k$ large. 
We can now state the main result.
  \begin{theorem}[Main result] We have, for all $f \in C^0(M)$,
$$ \mathcal{H}^{n-1} \left\{x:f(x) =0\right\} \gtrsim_{(M,g)}  \left( \frac{ \|f\|_{L^1}}{\|f\|_{L^{\infty}}} \right)^{2 - \frac{1}{n}}    \frac{  \sqrt{\lambda}}{(\log{\lambda})^{n/2}}.$$
  \end{theorem}
\vspace{10pt}
{\noindent \textit{Comments.} }Several comments and remarks are in order.\\

\textbf{1.} The scaling in $\|f\|_{L^1}/\|f\|_{L^{\infty}}$ is most likely not optimal,
we have some comments on this after Theorem 2. The trigonometric function $f(x) = \cos{(k x_1)}$ on the flat torus $\mathbb{T}^d$ shows that the result is sharp up to the logarithmic factor since $\|f\|_{L^1}/\|f\|_{L^{\infty}} \sim 1$, $\mathcal{H}^{n-1}\left\{x:f(x)=0\right\} \sim k$ and $\lambda \sim k^2$.
This extends to 'flat functions' $c^{-1} \leq \|f\|_{L^{1}}, \|f\|_{L^{\infty}} \leq c$ on general manifolds: if $\left\langle f, \phi_n\right\rangle = 0$ for all eigenfunctions $\phi_n$ with eigenvalue $\lambda_n \leq \lambda$
$$ \mathcal{H}^{n-1} \left\{x:f(x) =0\right\} \gtrsim_{(M,g),c}   \frac{  \sqrt{\lambda}}{(\log{\lambda})^{n/2}}.$$
It seems likely that the logarithm is an artifact of the proof and can be removed after which the result would be sharp.\\

\textbf{2.} We always work with functions that are continuous since $\left\{x:f(x) = 0\right\}$ might be empty otherwise. However, it is not very important for the proof and we could relax the condition to $f\in L^{\infty}(M)$
provided that we reinterpret
$$ \mathcal{H}^{n-1}\left\{x:f(x)=0\right\} := \liminf_{t \rightarrow 0}{ \mathcal{H}^{n-1} \left\{x: \left[e^{t \Delta} f\right](x) = 0 \right\}},$$
where $e^{t\Delta}$ denotes the heat semigroup.\\

\textbf{3.} It is not difficult to see that for any set of eigenfunctions $\left\{ \phi_j: j \in J\right\}$ of eigenfunctions,
any linear combination has to vanish somewhere (because of orthogonality to constants). One particular consequence is that we are able to determine bounds on the size of the set where they vanish. 

\begin{corollary} Let $J \subset \mathbb{N}$. Then, for all $a_j \in \mathbb{R}$, and all $\varepsilon>0$
$$ \mathcal{H}^{n-1} \left\{x \in M: \sum_{j \in J}{a_j \phi_j} = 0 \right\} \gtrsim_{(M,g), \varepsilon}  \left( \frac{\min_{j \in J} \lambda_j^{\frac{n}{4n-2}}}{ \max_{j \in J}{ \|\phi_j\|_{L^{\infty}}^{2}} }\frac{\sum_{j \in J}{|a_j|^{2}}}{\left(\sum_{j \in J}{|a_j|}\right)^{2}}\right)^{2-\frac{1}{n} + \varepsilon}.$$
\end{corollary}
This is certainly very far from a sharp result.
It is tempting to believe that the highest oscillatory term dominates the entire linear expression and that, for example, 
$$ \mathcal{H}^{n-1}\left\{x: \phi_{\lambda}(x) = \phi_{\mu}(x)\right\} \gtrsim \max\left( \sqrt{\lambda}, \sqrt{\mu} \right)$$
uniformly in $\mu, \lambda$. Corollary 1 implies only the much weaker result
$$ \mathcal{H}^{n-1}\left\{x: \phi_{\lambda}(x) = \phi_{\mu}(x)\right\} \gtrsim_{\varepsilon} \frac{\min\left( \lambda, \mu \right)^{1/2 - \varepsilon}}{\max\left( \lambda, \mu \right)^{n - \frac{3}{2} + \frac{1}{2n} + \varepsilon}}.$$

\subsection{A refinement for regular zero sets.}
The proof of Theorem 1 has a combinatorial ingredient coming from the fact that the set $\left\{x: f(x) = 0\right\}$ can be rather complicated. 
 If we can assume that it is not too irregular, then we can avoid combinatorial reasoning and stronger results are possible.
$\left\{x: f(x) = 0\right\}$ being not 'too irregular' is defined via volume-expansion: we ask that the volume of a small $\varepsilon-$neighborhood of $\left\{x:f(x) = 0\right\}$ grows like 
$ \sim \varepsilon \mathcal{H}^{n-1}\left\{x:f(x) = 0\right\}$ for small $\varepsilon$.

\begin{center}
\begin{figure}[h!]
\begin{tikzpicture}[scale=4]
\draw [ultra thick] (0,0) -- (1,0) -- (1,1) -- (0,1) -- (0,0);
\draw [ thick] (0, 0.3) to[out=0, in=120] (0.4, 0);
\draw [ thick] (0.4, 1) to[out=330, in=180] (1, 0.3);
\draw [ thick] (0,1) -- (1,0);
\draw [ thick] (0,0.6) -- (0.6,0);
\draw [ thick] (0,0) -- (1,1);

\draw [ultra thick] (1.6,0) -- (2.6,0) -- (2.6,1) -- (1.6,1) -- (1.6,0);
\draw [ thick] (1.7, 0.2) circle (0.01cm);
\draw [ thick] (1.8, 0.4) circle (0.01cm);
\draw [ thick] (2, 0.6) circle (0.01cm);
\draw [ thick] (2.1, 0.5) circle (0.01cm);
\draw [ thick] (2.4, 0.2) circle (0.01cm);
\draw [ thick] (2.3, 0.9) circle (0.01cm);
\draw [ thick] (1.8, 0.6) circle (0.01cm);
\draw [ thick] (2.0, 0.2) circle (0.01cm);
\draw [ thick] (2.2, 0.5) circle (0.01cm);
\end{tikzpicture}
\caption{Zero sets with slow (left) and fast (right) volume-expansion.}
\end{figure}
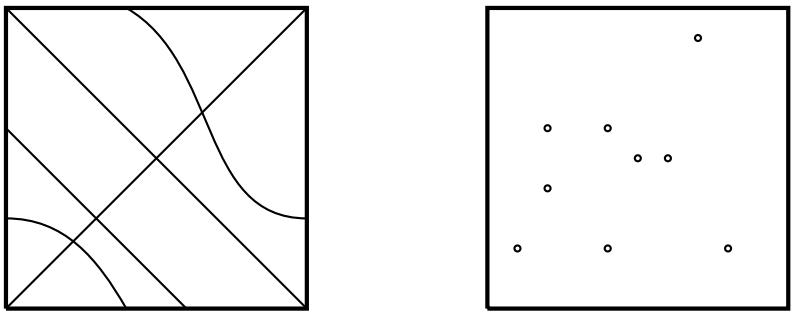
\end{center}
The condition is easily satisfied if the nodal set has large parts that are relatively flat or highly twisted but confined to a small region of space; it is violated if the nodal
set breaks up in a great number of small pieces that live below the wavelength.
\begin{theorem}
If the $\varepsilon-$neighborhood of $\left\{x: f(x) = 0\right\}$ satisfies
$$ \mathcal{H}^n\left( \left\{x: f(x) = 0\right\} + B_{\varepsilon}\right) \leq c \varepsilon\mathcal{H}^{n-1}\left\{x: f(x) = 0\right\}$$
for some $c>0$ and for all
$$0 \leq \varepsilon \lesssim \frac{1}{\sqrt{\lambda}} \log{\frac{\|f\|_{L^{\infty}}}{\|f\|_{L^1}}},$$
then
$$ \mathcal{H}^{n-1}\left\{x: f(x) = 0\right\} \gtrsim_{(M,g),c}  \frac{\|f\|_{L^1}}{\|f\|_{L^{\infty}}} \left(1 + \log{\frac{\|f\|_{L^{\infty}}}{\|f\|_{L^1}}} \right)^{-1} \sqrt{\lambda}.$$
\end{theorem}

The proof suggests that the logarithm is an artifact of the proof and might be removable. This result (without the logarithm) has such a simple form that we are tempted
to believe that it might point towards what could be true. A somewhat supporting heuristic is the following: let us consider, using $\delta_x$ to denote the Dirac measure in $x$, and for $0 \leq r \ll 1$,
$$ f = - n + e^{t\Delta} \sum_{k=1}^{n}{\delta_{x_k}}  \qquad \mbox{for}~t = r^2 n^{-2/d}.$$
Generically, we would expect to be able to pick $n$ points in such a way that $f$ is orthogonal to the first $\sim n$ Laplacian eigenfunctions (and one would not expect to
be able to do better than that, see \cite{steinm}). A simple computation
shows 
$$ \mathcal{H}^{n-1}\left\{x: f(x) = 0\right\}  \sim r^{d-1} n^{\frac{1}{d}}, \quad \|f\|_{L^1} \sim n,~\quad \mbox{and} \quad \|f\|_{L^{\infty}} \sim r^{-d} n.$$
Moreover, Weyl's law yields the scaling $\lambda \sim n^{2/d}$. This would suggest
$$ \mathcal{H}^{n-1}\left\{x: f(x) = 0\right\} \sim \left(\frac{\|f\|_{L^1}}{\|f\|_{L^{\infty}}}\right)^{\frac{d-1}{d}}\sqrt{\lambda}.$$
By letting $r \rightarrow 0$, the heuristic suggests that the bound cannot purely depend on $\lambda$ but will have to involve the function in some way. The example is not rigorous
since we do not currently know whether there exist $\sim n$ points that are well spaced out with the property that smoothed Dirac measures located in these points are orthogonal to the first $\sim n$ eigenfunctions.

\section{Proofs}
\subsection{Proof of Theorem 1}
The proof can be roughly summarized as follows: we decompose the manifold into many little boxes and, depending on the size of $\|f\|_{L^1}$ and $\|f\|_{L^{\infty}}$, can guarantee to find
at least a certain number of boxes with a nontrivial $L^1-$mass. Since $f$ has most of its $L^2-$energy at high frequencies, it decays rapidly under the heat flow. Interpreting the heat flow
as a local process arising from convolution with the heat kernel implies that there are a quite a few boxes with a lot of $L^1-$mass that becomes small when taking local averages. This means that
there has to be both positive and negative $L^1-$mass at a local scale for many different small boxes, which forces vanishing on a large set.\\

As for notation, we shall use $A \lesssim B$ and $A \sim B$ to hide the existence of universal constants (that may all depend on each other). Explicit numerical constants, if needed, are
chosen for simplicity of exposition and never designed to be sharp.

\begin{proof}[Proof of Theorem 1.] We assume w.l.o.g. that $\|f\|_{L^{\infty}} = 1$ and that $\vol(M) = 1$ (mainly to simplify book-keeping and being able to discard several unimportant constants). We will denote the heat evolution of $f$ by
$$e^{t\Delta}f = \sum_{k=0}^{\infty}{ e^{-\lambda_k t} \left\langle f, \phi_k\right\rangle \phi_k}.$$
Using the Cauchy-Schwarz inequality and $L^2-$orthogonality, we obtain
\begin{align*}
\| e^{t\Delta} f\|_{L^1} \leq \| e^{t\Delta} f\|_{L^2} = \left(\left\| \sum_{\lambda_k \leq \lambda}{ e^{-\lambda_k t} \left\langle f, \phi_k \right\rangle \phi_k }\right\|^2_{L^2} + \left\| \sum_{\lambda_k \geq \lambda}{ e^{-\lambda_k t} \left\langle f, \phi_k \right\rangle \phi_k }\right\|^2_{L^2}\right)^{1/2}.
\end{align*}
The elementary inequality $(a^2 + b^2)^{1/2} \leq a + b$ for $a,b > 0$, the definition of $\lambda$ and $\|f\|_{L^{\infty}}=1$ then imply
\begin{align*}
\| e^{t\Delta} f\|_{L^1} &\leq \left\| \sum_{\lambda_k \leq \lambda}{ \left\langle f, \phi_k \right\rangle \phi_k }\right\|_{L^2} + e^{-\lambda t} \left\| \sum_{\lambda_k \geq \lambda}{ \left\langle f, \phi_k \right\rangle \phi_k} \right\|_{L^2}\\
& \leq c\|f\|_{L^1} +  e^{-\lambda t}  \|f\|_{L^2}\\
&\leq c\|f\|_{L^1}+  e^{-\lambda t}   \|f\|_{L^1}^{1/2}. 
\end{align*}
This means that as soon as $t$ is so large that
$$ e^{-\lambda t} \ll \|f\|_{L^1}^{1/2} \qquad \mbox{we have} \qquad \| e^{t\Delta} f\|_{L^1}  \sim c \| f\|_{L^1}.$$
This fixes the natural time-scale 
$$ t  \sim \frac{1}{\lambda} \log{\left(\frac{e}{\|f\|_{L^1}}\right)}.$$
The implicit constant is chosen so large that $\| e^{t\Delta} f\|_{L^1}  \leq (c_n^{}/10000)  \| f\|_{L^1},$ where $0<c_n \ll 1 $ is a constant depending only on the dimension that will be introduced later (and, in
particular, $c=c_n/10000$ in the definition of $\lambda$ will work).
The heat equation can be written as a convolution
$$ e^{t\Delta}f(x) = \int_{M}{ p_t(x,y) f(y) dy},$$
where the heat kernel satisfies a Gaussian bound of the type
$$  p_t(x,y) \leq \frac{c_1}{t^{n/2}} \exp\left( -c_2 \frac{d(x,y)^2}{t}\right),$$
where $c_1, c_2 >0$ are constant depending only on $(M,g)$. This means that at time $t$, the heat evolution $e^{t \Delta} f$ can be thought of as a local averaging at scale $\sim \sqrt{t}$.
The second part of the argument is more combinatorial. We introduce the spatial scale $\delta = \sqrt{t}$ and partition $(M,g)$ into non-overlapping
cube-like sets $Q_i$ all of which have the same measure $\delta^n$ and diameter $\lesssim \delta$. 
Clearly, one would expect that massive global decay requires decay on many of the $Q_i$. We will now argue that there is a substantial number of sets $Q_i$ on which the $L^1-$norm is not too
small and where we observe a nontrivial amount of decay.\\

For a parameter $0 < c_n \ll 1/1000$ (same as above, depending only on the dimension and introduced below), we define the sets
\begin{align*}
A &=  \left\{Q_i : \frac{1}{|Q_i|}\int_{Q_i}{|f| dx}\ \leq  \frac{1}{2} \|f\|_{L^1} ~\mbox{and}~ \int_{Q_i}{|e^{t \Delta} f| dx} \leq \frac{c_n}{100}\int_{Q_i}{|f| dx}\right\} \\
B &= \left\{ Q_i : \frac{1}{|Q_i|}\int_{Q_i}{|f| dx}\ >  \frac{1}{2}\|f\|_{L^1} ~\mbox{and}~  \int_{Q_i}{|e^{t \Delta} f| dx} \leq  \frac{c_n}{100}\int_{Q_i}{|f| dx}\right\}\\
C &= \left\{ Q_i:  \int_{Q_i}{|e^{t \Delta} f| dx} \geq \frac{c_n}{100} \int_{Q_i}{|f| dx}\right\}
\end{align*}
and want to show that the set $B$ has to be large. We start by showing that not a lot of $L^1-$mass can be in the set $C$. Clearly,
$$\| e^{t\Delta} f\|_{L^1}  \leq \frac{c_n^{}}{10000}  \| f\|_{L^1} \qquad \mbox{implies} \qquad \sum_{Q_i \in C}{\int_{Q_i}{ |f| dx}} \leq \frac{1}{100} \|f\|_{L^1}.$$
As a consequence of $|f| \leq 1$ and $|Q_i| = \delta^n$
$$ \left(1 - \frac{1}{100} \right) \|f\|_{L^1} \leq \sum_{Q_i \in A \cup B}{\int_{Q_i}{ |f| dx}} \leq \frac{1}{2} \|f\|_{L^1} \delta^n (\# A) + \delta^n (\# B)$$
and thus, since $(\#A)\delta^{-n} \leq 1$,
 $$\#B \gtrsim \delta^{-n} \|f\|_{L^1}.$$
We will now introduce the geometric constant $c_n$ (depending on the dimension) as follows: given fixed $(M,g)$, we can achieve a sufficiently fine
partition into 'almost' cubes. Assuming minimal regularity on the almost-cubes (which can be achieved to an arbitrary degree of accuracy for $\delta$
sufficiently small), we obtain for all $Q_i$, using $\delta_x$ to denote the Dirac measure in $x$,
$$ \inf_{x \in Q_i} \int_{Q_i}{ e^{t\Delta}\delta_x} \geq c_n$$
with $c_n$ depending only on the dimension.
\begin{center}
\begin{figure}[h!]
\begin{tikzpicture}[scale=2]
\draw [ultra thick] (0, 0) to[out=0, in=180] (1, 0.2) to[out=80,in=270]  (1.1, 1) to[out=180, in=0] (0.1, 0.9) to[out=270, in=90] (0,0);
\filldraw [ultra thick] (0,0) circle (0.03cm);
\node at (-0.1, 0) {$x$};
\draw [ultra thick] (0,0) circle (0.3cm);
\draw [ultra thick] (0,0) circle (0.4cm);
\draw [ thick]  (0,0) circle (0.5cm);
\draw [ thick]  (0,0) circle (0.6cm);
\draw   (0,0) circle (0.7cm);
\draw (0,0) circle (0.8cm);
\draw [dashed] (0,0) circle (0.9cm);
\draw [dashed] (0,0) circle (1cm);
\end{tikzpicture}
\caption{The constant $c_n$: for a cube-like decomposition, the heat of a point mass in $x \in Q_i$ has to localize
at least $c_n$ of its mass inside the cube.}
\end{figure}
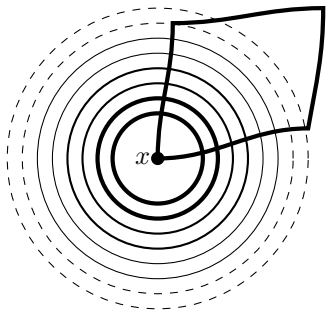
\end{center}

It seems likely that one could actually achieve, up to relatively small errors decreasing with $\delta$,
$$ c_n \sim \frac{1}{(2\pi)^{n/2}}\int_{[0,1]^n}{ e^{-\frac{\|x\|^2}{2}} dx},$$
however, this will not be of importance. (Likewise, the cubical shape is not important and cubes could be replaced by sets with bounded diameter.)
 If $c_n \geq 1/1000$, which is likely in low dimensions, then we make it artificially smaller and set $c_n =1/10000$ (its precise value is not important as long
as it is sufficiently small).
 One particularly useful consequence is
that it allows to conclude, for any nonnegative $g \geq 0$, that
$$ \int_{Q_i}{e^{t\Delta}(g \chi_{Q_i}) dx} \geq c_n \int_{Q_i}{g dx},$$
which means that we can control the loss of mass to the outside.
Let us now consider such a $Q_i$ for some $Q_i \in B$ and assume w.l.o.g. 
$$ \int_{Q_i}{\max(f,0)dx} \geq \int_{Q_i}{ |\min(f,0)|dx}.$$ 
Then, since $Q_i \in B$,
$$  \int_{Q_i}{\max(f,0)dx}  \geq \frac{1}{4} \|f\|_{L^1} \delta^n \qquad \mbox{and} \qquad \int_{Q_i}{|e^{t\Delta} f|dx}  \leq \frac{c_n}{50}\int_{Q_i}{\max(f,0)dx}.$$
Writing the heat evolution as a convolution, we observe that for the averaging to be small, we require cancellation. We will now argue that
there has to be at least some negative mass in a $r-$neighborhood of $Q_i$, for $r$ to be determined,
$$ \int_{d(x,Q_i) \lesssim r}{|\min(0, f(x))|dx}  \geq   \frac{c_n}{10000}  \|f\|_{L^1} \delta^n.$$
We shall now determine an upper bound on how large $r$ needs to be for this inequality to be true. Let us assume the statement is false and
$$ \int_{d(x,Q_i) \lesssim r}{|\min(0, f(x))|dx}  \leq \frac{c_n}{10000} \|f\|_{L^1} \delta^n.$$
Our goal is to find $r$ that makes this inequality incorrect: to do so, we now try to bound the amount of decay happening within $Q_i$ in the most favorable way. This is the following: all the negative $L^1-$mass in $d(x, Q_i) \lesssim r$ flows to $Q_i$ and
cancels with positive mass, the positive mass in $Q_i$ diffuses to the outside as much as possible,
there is no other positive mass anywhere on the manifold and the rest of the $L^1-$mass is negative
and localized at distance $\sim r$ from $Q_i$.
We argue
$$ e^{t\Delta}f \geq e^{t\Delta}\max(f,0) \chi_{Q_i} + e^{t\Delta}\min(f,0) $$
and observe that
\begin{align*}
 \int_{Q_i}{e^{t\Delta}\max(f,0) \chi_{Q_i} dx} &\geq   \left(   \inf_{x \in Q_i} \int_{Q_i}{ e^{t\Delta}\delta_x} \right)   \int_{Q_i}{\max(f,0)dx}  \\
&\geq c_n \int_{Q_i}{\max(f,0)dx} \geq \frac{c_n}{4} \|f\|_{L^1} \delta^n.
\end{align*}
At the same time, 
\begin{align*}
 \int_{Q_i}{e^{t\Delta}\min(f,0)\chi_{d(x,Q_i) \lesssim r} dx} \leq  \int_{d(x,Q_i) \lesssim r}{|\min(0, f(x))|dx}  \leq \frac{c_n}{10000} \|f\|_{L^1} \delta^n
\end{align*}
and, using the decay of the heat kernel,
\begin{align*}
  \int_{Q_i}{e^{t\Delta}\min(f,0)\chi_{d(x,Q_i) \gtrsim r} dx} &\geq  \frac{c_1}{t^{n/2}} \exp{\left( -c_2 \frac{r^2}{t} \right)}\int_{Q_i}{\min(f,0)\chi_{d(x,Q_i) \gtrsim r} dx} \\
&\geq -\frac{c_1}{t^{n/2}} \exp{\left( -c_2 \frac{r^2}{t} \right)}\|f\|_{L^1}
\end{align*}
Altogether, we obtain the estimate
\begin{align*}
\int_{Q_i}[e^{t \Delta}f](x)dx \geq \frac{c_n}{5} \|f\|_{L^1} \delta^n - \frac{c_1}{t^{n/2}}\exp\left(-c_2 \frac{r^2}{t}\right)\|f\|_{L^1}.\end{align*}
At the same time, since $Q_i \in B$, we know that there is decay by at least a factor of $c_n/50$ and thus
$$ \frac{c_1}{t^{n/2}}\exp\left(-c_2 \frac{r^2}{t}\right)\|f\|_{L^1} \gtrsim c_n \|f\|_{L^1}\delta^n \gtrsim_{(M,g)} \|f\|_{L^1}\delta^n .$$
Since $t^{n/2} = \delta^n$, we see that the appropriate scale for $r$ is given by 
$$ r \lesssim \delta \left(\log{ \frac{1}{\delta}}\right)^{1/2}.$$
The isoperimetric inequality implies, for $Q_i \in B$,
$$ \mathcal{H}^{n-1}\left\{x: d(x,Q_i) \lesssim \delta \left(\log{  \frac{1}{\delta} }  \right)^{1/2}~ \mbox{and} ~ f(x) = 0\right\} \gtrsim  \|f\|^{\frac{n-1}{n}}_{L^1} \delta^{n-1}.$$
We now want to find isolated cubes $Q_i$ that are at least $r-$separated and then apply this inequality. A simple volume comparison shows that we can find at least 
$$ \frac{(\#B) \delta^n}{ r^n} \gtrsim \frac{\#B}{  \left(\log{\frac{1}{\delta} }\right)^{n/2}} \gtrsim  \frac{\delta^{-n} \|f\|_{L^1} }{  \left(\log{  \frac{1}{\delta}  }\right)^{n/2}}$$
such cubes and this implies
$$ \mathcal{H}^{n-1}\left\{x: f(x) = 0\right\} \gtrsim  \|f\|^{\frac{n-1}{n}}_{L^1} \delta^{n-1}  \frac{\delta^{-n} \|f\|_{L^1} }{  \left(\log{   \frac{1}{\delta}  }\right)^{n/2}} = \frac{1}{\delta}\frac{1}{{  \left(\log{  \frac{1}{\delta} }\right)^{n/2}}} \|f\|_{L^1}^{2-\frac{1}{n}} .$$
The desired result then follows from recalling that
$$ \delta = \frac{1}{\sqrt{\lambda}} \left( \log \frac{e}{\|f\|_{L^1}}\right)^{1/2} \gtrsim \frac{1}{\sqrt{\lambda}}.$$
\end{proof}

\subsection{Proof of Theorem 2}
The argument is similar in spirit to that of the main result but much simpler: if the $\varepsilon-$neighborhoods of $\left\{x:f(x) = 0\right\}$ are small, then the normalization
$\|f\|_{L^{\infty}}=1$ implies that only a certain proportion of the $L^1-$mass can be close to $L^1-$mass of the opposite sign. Since, by the same reasoning as above, we have rapid
decay of the heat semigroup in time, i.e. for
$$ t  \sim \frac{1}{\lambda} \log{\left(\frac{e}{\|f\|_{L^1}}\right)}$$
with a suitable implicit constant, we get
$$ \| e^{t\Delta} f\|_{L^1} \leq \| e^{t\Delta} f\|_{L^2}  \leq \frac{\|f\|_{L^1}}{10000}+  e^{-\lambda t}   \|f\|_{L^1}^{1/2},$$
this can be used to force a contradiction if $\left\{x:f(x) = 0\right\}$ is too small. A suitable metaphor is that opening the windows to cool a house is only effective if the windows are not too small.
The new mathematical ingredient is the Davies-Gaffney estimate \cite{davies,gaffney,grig}:
if $A, B$ are measurable subsets of $M$, then for all $t>0$
$$ \int_{A}\int_{B} p_t(x,y) dg dg \leq \sqrt{|A||B|} \exp\left(-\frac{d(A,B)^2}{4t}\right),$$
where $d(A,B) = \inf_{a \in A, b\in B}{d(a,b)}.$ We use the version for $L^2-$functions (see \cite{grig2}) 
$$ |\left\langle e^{t\Delta} f_1, f_2\right\rangle | \leq \exp\left(-\frac{d(\mbox{supp}(f_1),\mbox{supp}(f_2))^2}{4t}\right) \|f_1\|_{L^2} \|f_2\|_{L^2},$$
where $\mbox{supp}$ denotes the support of the function.

\begin{proof}[Proof of Theorem 2]

Let us assume again w.l.o.g. that $\|f\|_{L^{\infty}} = 1$. We abbreviate $\Sigma = \left\{x: f(x) = 0\right\}$ and try to get an estimate on its size. 
The assumption on the regularity of the nodal set implies that
$$ \left\{x: d(x,Z) \leq \delta\right\} \leq c \delta \mathcal{H}^{n-1} (\Sigma).$$
Our assumption $\|f\|_{L^{\infty}}=1$ implies that whenever $c \delta  \mathcal{H}^{n-1} (\Sigma) \leq \|f\|_{L^1}/10,$ then
$$ \int_{  d(x,Z) \leq \delta }{|f|dx} \leq \int_{d(x,Z) \leq \delta  }{\|f\|_{L^{\infty}}dx} = c \delta  \mathcal{H}^{n-1} (\Sigma) \leq  \frac{\|f\|_{L^1}}{10},$$
which shows that 90\% of the $L^1-$mass is at least distance $\delta$ away from $\Sigma$.
We define four regions
\begin{align*}
 A &= \left\{x \in M: d(x, \Sigma) \geq \delta \wedge f(x) < 0\right\} \\
 B &= \left\{x \in M: d(x, \Sigma) \leq \delta \wedge f(x) < 0\right\} \\
 C &= \left\{x \in M: d(x, \Sigma) \leq \delta \wedge f(x) > 0\right\}\\
 D &= \left\{x \in M: d(x, \Sigma) \geq \delta \wedge f(x) > 0\right\}
\end{align*}
and will now 'manually' control the flow of the $L^1-$mass. Since $\lambda > 0$, we have almost orthogonality to constants 
$$ |\left\langle f, 1\right\rangle| \leq \frac{ \|f\|_{L^1}}{10000}$$
and thus
$$ \left| \int_{A \cup B}{|f| dx} - \int_{C \cup D}{|f|dx} \right| \leq \frac{1}{100}\|f\|_{L^1}.$$

\begin{center}
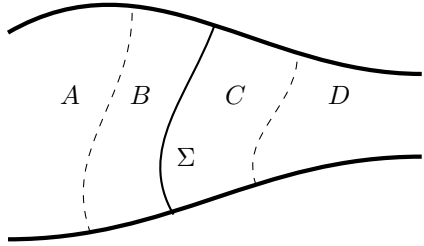
\begin{figure}[h!]
\begin{tikzpicture}[scale=5.5]
\draw [ultra thick] (0, 0) to[out=0, in=180] (1, 0.2);
\draw [ultra thick] (0, 0.5) to[out=30, in=180] (1, +0.4);
\draw [thick] (0.4, 0.06) to[out=120, in=250] (0.5, +0.52);
\node at (0.43, 0.2) {$\Sigma$};
\draw [dashed] (0.6, 0.13) to[out=110, in=270] (0.7, +0.45);
\draw [dashed] (0.2, 0.013) to[out=110, in=270] (0.3, +0.55);
\node at (0.15, 0.35) {$A$};
\node at (0.32, 0.35) {$B$};
\node at (0.55, 0.35) {$C$};
\node at (0.8, 0.35) {$D$};
\end{tikzpicture}
\caption{The local picture: the separating set $\Sigma$ and regions $\delta-$close and further away from it. It is difficult for mass from $A$ to flow to $C \cup D$ if $B$ is wide.}
\end{figure}
\end{center}

Our goal is an estimate
$$ \int_{C \cup D}{e^{t \Delta} f dx} \qquad \mbox{from below.}$$
We will obtain a bound from below that gets smaller as $\mathcal{H}^{n-1}(\Sigma)$ gets bigger. Simultaneously, because $f$ lives at high-frequencies,
we know that there is a fixed amount of decay. That will then imply a lower bound on $\mathcal{H}^{n-1}(\Sigma)$. 
We approach this using linearity
$$ e^{t\Delta} f = e^{t\Delta} f \chi_{A} + e^{t\Delta} f \chi_{B}  + e^{t\Delta} f \chi_{C}  + e^{t\Delta} f \chi_{D}$$
and analyzing the flows of each of the four solutions separately. The worst case is when all the positive $L^1-$mass
from $C$ wanders into $A \cup B$ and the negative $L^1-$mass from $B$ wanders into $C \cup D$. It is also
harmful for positive $L^1-$mass from $D$ to wander into $A \cup B$ and, conversely, for negative $L^1-$mass
from $A$ to wander into $C \cup D$, however, this is more difficult to do since they have cross a certain distance
first. More formally, the first two steps translate into the trivial estimates
$$ \int_{C \cup D}{ e^{t\Delta} (f \chi_{B} )dx} \geq \int_{B}{ fdx} \qquad \mbox{and} \qquad  \int_{C \cup D}{ e^{t\Delta} (f \chi_{C} )dx} \geq 0.$$
For the other two steps, we use the Davies-Gaffney inequality to bound
$$\left| \left\langle e^{t\Delta} (f \chi_{A}), \chi_{C \cup D}\right\rangle \right| \leq e^{-\frac{\delta^2}{4t}} \|f\chi_A\|_{L^2}  |C \cup D|^{1/2}  \leq   e^{-\frac{\delta^2}{4t}} \|f\chi_A\|_{L^2} $$
and thus
$$ \int_{C \cup D}{e^{t\Delta} (f \chi_{A} ) dx} \geq  -   e^{-\frac{\delta^2}{4t}} \|f\chi_A\|_{L^2} .$$
Likewise, the very same reasoning also yields
$$ \int_{A \cup B}{e^{t\Delta} (f \chi_{D} ) dx} \leq    e^{-\frac{\delta^2}{4t}}\|f\chi_D\|_{L^2} .$$
Moreover, the assumption of neighborhoods of $\Sigma$ not growing too rapidly combined with $\|f\|_{L^{\infty}} = 1$ implies
$$ \int_{B}{|f|dx} \leq |B| \leq |B \cup C| \leq  c \delta  \mathcal{H}^{n-1} (\Sigma).$$
Now we can collect all these arguments. We have
\begin{align*}
\int_{C \cup D}{e^{t \Delta} f dx} &\geq \int_{ D}{f dx} -  \int_{B}{ |f|  dx} -  e^{-\frac{\delta^2}{4t}}  \|f\chi_A\|_{L^2} -  e^{-\frac{\delta^2}{4t}} \|f\chi_D\|_{L^2} \\
&\geq  \int_{ D}{f dx}  - 2 e^{-\frac{\delta^2}{4t}} \|f\|_{L^1}^{1/2} - c \delta  \mathcal{H}^{n-1} (\Sigma).
\end{align*}
Clearly, by our considerations above,
$$ \int_{C \cup D}{e^{t \Delta} f dx} \leq \int_{M}{|e^{t \Delta} f| dx} \leq  \frac{\|f\|_{L^1}}{10000}+  e^{-\lambda t}   \|f\|_{L^1}^{1/2}.$$
Combining the last two inequalities yields
\begin{align*}
c \delta \mathcal{H}^{n-1}(\Sigma) \geq \int_{D}{f dx} - \frac{\|f\|_{L^1}}{10000} - \left( 2 e^{-\frac{\delta^2}{4t}} + e^{-\lambda t} \right) \|f\|_{L^1}^{1/2}.
\end{align*}
The volume-expansion bound, almost orthogonality to constants and $\|f\|_{L^{\infty}}=1$ imply 
$$\int_{D}{f dx} \geq  \frac{4}{10} \|f\|_{L^1} - \int_{C}{f dx} \geq \frac{4}{10} \|f\|_{L^1} - c \delta \mathcal{H}^{n-1}(\Sigma) $$
and thus 
$$ 2c \delta \mathcal{H}^{n-1}(\Sigma) \geq \frac{4}{10} \|f\|_{L^1} - \left( 2 e^{-\frac{\delta^2}{4t}} + e^{-\lambda t} \right) \|f\|_{L^1}^{1/2} $$

This allows us now to pick suitable parameters. Clearly, we want to choose $t$ and $\delta$ (depending on $t$) such that the negative
quantities on the right-hand side are $\lesssim \|f\|_{L^1}$ to obtain a nonnegative right-hand side. Formally, we need
$$  e^{-\lambda t}  \leq \frac{ \|f\|^{1/2}_{L^1}}{1000} \qquad \mbox{forcing} \qquad t \gtrsim  \frac{1}{\lambda} \log{\left(\frac{1}{\|f\|_{L^1}}\right)}$$
and then
$$ 2 e^{-\frac{\delta^2}{4t}} \leq \frac{ \|f\|^{1/2}_{L^1}}{1000}  \qquad \mbox{forcing} \qquad \delta \gtrsim \sqrt{t} \left[ \log{\left(\frac{1}{\|f\|_{L^1}}\right)}\right]^{1/2}.$$
Picking these values of $t$ and $\delta$ then implies
$$ \delta \sim \frac{1}{\sqrt{\lambda}} \log{\left(\frac{1}{\|f\|_{L^1}}\right)}$$
and this yields
$$ \mathcal{H}^{n-1}(\Sigma) \gtrsim \frac{1}{\delta} \|f\|_{L^1}.$$
\end{proof}

\subsection{Proof of Corollary 1}

\begin{proof} The argument is completely straightforward. 
Clearly,
$$ \xi \left(\sum_{j \in J}{a_j \phi_j}\right) = \min_{j \in J}{\lambda_j}.$$
$L^2-$orthogonality implies
$$ \sum_{j \in J}{a_j^2} = \left\| \sum_{j \in J}{a_j \phi_j} \right\|_{L^2}^2 \leq \left\| \sum_{j \in J}{a_j \phi_j} \right\|_{L^1}  \left\| \sum_{j \in J}{a_j \phi_j} \right\|_{L^{\infty}}.$$
The triangle inequality yields
$$  \left\| \sum_{j \in J}{a_j \phi_j} \right\|_{L^{\infty}} \leq  \left(\max_{j \in J}{\|\phi_j\|_{L^{\infty}}}\right) \sum_{j \in J}{|a_j |}.$$
Altogether, this implies
$$ \frac{\|f\|_{L^1}}{\|f\|_{L^{\infty}}} \gtrsim \frac{1}{ \left(\max_{j \in J}{\|\phi_j\|_{L^{\infty}}}\right) ^2}  \frac{\sum_{j \in J}{a_j^2}}{\left(\sum_{j \in J}{|a_j|}\right)^2}.$$
\end{proof}

\end{document}